\documentclass[11pt,leqno]{amsart}
  
\usepackage{latexsym,enumitem}
\usepackage{amssymb}
\usepackage[cp850]{inputenc}
\usepackage{epsfig}
\usepackage{psfrag}
\usepackage{amsthm}
\usepackage{amscd}
\usepackage{amsmath}
\usepackage{amsfonts}
\usepackage{graphics,caption}
\usepackage[all]{xy}
\usepackage{etoolbox}
\patchcmd{\quote}{\rightmargin}{\leftmargin 2em \rightmargin}{}{}
 \usepackage{mathtools}

\usepackage[english]{babel}
\usepackage{comment}
\usepackage{color}
\usepackage{hyperref}

\let\phi\varphi

\DeclareMathOperator{\eqdef}{\coloneqq} %definition
\let\epsilon\varepsilon
\let\subset\subseteq
\newcommand{\be}{\begin{equation*}}
 \newcommand{\ee}{\end{equation*}}
\newcommand{\bpf}{\begin{dimo}}
\newcommand{\epf}{\end{dimo}}
\newcommand{\bdefi}{\begin{defin}}
\newcommand{\edefi}{\end{defin}}
\newcommand{\bthm}{\begin{thm}}
\newcommand{\ethm}{\end{thm}}
\newcommand{\blem}{\begin{lem}}
\newcommand{\elem}{\end{lem}}
\newcommand{\bcor}{\begin{cor}}
\newcommand{\ecor}{\end{cor}}

\newcommand{\bprop}{\begin{prop}}
\newcommand{\eprop}{\end{prop}}
\newcommand{\bese}{\begin{ese}}
\newcommand{\eese}{\end{ese}}
\newcommand{\brem}{\begin{rem}}
\newcommand{\erem}{\end{rem}}
\newcommand{\bpfc}{\begin{dimoclaim}}
\newcommand{\epfc}{\end{dimoclaim}}

\newcommand{\rar}{\rightarrow} %arrow right
\newcommand{\Orar}[1]{\xrightarrow{{}_{#1}}} %long arrow right with overset
\newcommand{\imm}{\looparrowright}

\newcommand{\abs}[1]{\left\lvert#1\right\rvert}						%mod
\newcommand{\set}[1]{\left\{#1\right\}}					%set, curly brackets
\DeclareMathOperator{\diam}{diam}					%Diameter
\newcommand{\quotient}[2]{\left.\raisebox{.1em}{$#1\!$}\middle/\raisebox{-.1em}{$#2$}\right.}

\DeclareMathOperator{\emp}{\varnothing} %Known Sets, Fields and so on...
\DeclareMathOperator{\N}{\mathbb N}			%Naturals (with 0)
\DeclareMathOperator{\R}{\mathbb R}			%Reals
\DeclareMathOperator{\C}{\mathbb C}			%Complexes
\DeclareMathOperator{\Ha}{\mathbb H}			%Hamiltonian Quaternions
\newcommand{\M}{\mathscr M}

\newcommand{\bH}{\mathbb H^3}
 \newcommand{\hyp}[1]{\quotient{\bH}{#1}}
  
\newcommand{\subgroup}{\leqslant}

\newenvironment{quot}
{
	\vspace{-0.2cm}
	\vspace{0.2cm}
}

\theoremstyle{definition}
\newtheorem{d1}{Definition}[section] %fa partire tutto da [section]

\newenvironment{defin}
{
	\begin{quot}
		\begin{d1}
		}
		{\end{d1}
	\end{quot}

}

\theoremstyle{definition}
\newtheorem{r1}[d1]{Remark}%[section]

\newenvironment{rem}
{
	\begin{quot}
		\begin{r1}
		}
		{\end{r1}
	\end{quot}
}

\theoremstyle{definition}
\newtheorem{e1}[d1]{Exercise}%[section]

\theoremstyle{definition}
\newtheorem{ese1}[d1]{Example}

\newenvironment{ese}
{
	\begin{quot}
		\begin{ese1}
	}
	{	
		\end{ese1}
	\end{quot}
}

\theoremstyle{definition}

\theoremstyle{definition}
\newtheorem{f2}[d1]{Fact}

\theoremstyle{definition}

\theoremstyle{definition}

\theoremstyle{definition}
\newtheorem{t1}[d1]{Theorem}%[section]

\newenvironment{thm}
{
	\begin{quot}
		\begin{t1}}
		{\end{t1}
	\end{quot}
}

\theoremstyle{definition}
\newtheorem*{T1*}{Theorem}%[section]

\newenvironment{teor*}
{
	\begin{quot}
		\begin{T1*}}
		{\end{T1*}
	\end{quot}
}

\newenvironment{dimo}
{\begin{proof}[Proof]
	}
	{\end{proof}}

\newenvironment{dimoclaim}{\emph{Proof of Claim:}\;}{\hfill$\square$}

	\theoremstyle{definition}
	\newtheorem{l1}[d1]{Lemma}%[section]
	
	\newenvironment{lem}
	{
		\begin{quot}
			\begin{l1}}
			{\end{l1}
		\end{quot}
	}
	\theoremstyle{definition}
	\newtheorem{p1}[d1]{Proposition}%[section]
	
	\newenvironment{prop}
	{
		\begin{quot}
			\begin{p1}}
			{\end{p1}
		\end{quot}
	}
	
	\theoremstyle{definition}
	\newtheorem{c1}[d1]{Corollary}%[section]
	
	\newenvironment{cor}
	{
		\begin{quot}
			\begin{c1}}
			{\end{c1}
		\end{quot}
	}

		\renewenvironment{abstract}
	{\list{}{\rightmargin\leftmargin}%
		\item[\textbf{Abstract:}]\relax}
	{\endlist}

\newtheorem{innercustomthm}{Theorem}
\newenvironment{customthm}[1]
  {\renewcommand\theinnercustomthm{#1}\innercustomthm}
  {\endinnercustomthm}

\newtheorem{innercustomprop}{Proposition}
\newenvironment*{customprop}[1]
  {\renewcommand\theinnercustomprop{#1}\innercustomprop}
  {\endinnercustomprop}

 \newtheorem*{Theorem*}{Theorem}
 \newtheorem*{Proposition*}{Proposition}
 \newtheorem*{Lemma*}{Lemma}	
 \newtheorem*{Definition*}{Definition}

 \usepackage{float}
 \usepackage{graphicx}
 \usepackage{color}
 \usepackage{transparent}
 \usepackage{cite}
 \usepackage{mathtools}
 
\usepackage{xcolor} 
 \usepackage{tcolorbox}
\newtcbox{\hl}[1][red]{on line, arc=7pt,colback=#1!10!white,colframe=#1!50!black,
  before upper={\rule[-3pt]{0pt}{10pt}},boxrule=1pt, boxsep=0pt,left=6pt,
  right=6pt,top=2pt,bottom=2pt}

 \usepackage[numbers]{natbib}

  \renewcommand{\S}{\mathbb S^3}

   \renewcommand{\tilde}{\widetilde}
   \renewcommand{\hat}{\widehat }

\renewcommand{\M}{\mathcal M}

\renewcommand{\P}{\mathcal P}

	\renewenvironment{abstract}
	{\list{}{\rightmargin\leftmargin}%
		\item[\textbf{Abstract:}]\relax}
	{\endlist}

	\makeatletter
\def\subsection{\@startsection{subsection}{2}%
  \z@{.5\linespacing\@plus.7\linespacing}{.5\linespacing}%
  {\normalfont\bfseries}}
\makeatother

\makeatletter
\def\section{\@startsection{section}{2}%
  \z@{.5\linespacing\@plus.7\linespacing}{1\linespacing}%
  {\normalfont\bfseries}}
\makeatother

\makeatletter
\def\paragraph{\@startsection{paragraph}{4}%
  \z@\z@{-\fontdimen2\font}%
  {\normalfont\bfseries}}
\makeatother

\begin{document}
	 
\title{Hyperbolic Limits of Cantor Set complements in the Sphere}
\author{Tommaso Cremaschi and Franco Vargas Pallete}
\maketitle
\begin{abstract} Let $M$ be a hyperbolic 3-manifold with no rank two cusps admitting an embedding in $\S$. Then, if $M$ admits an exhaustion by $\pi_1$-injective sub-manifolds there exists Cantor sets $C_n\subset \S$ such that $N_n=\S\setminus C_n$ is hyperbolic and $N_n\rar M$ geometrically.
\end{abstract}

	\begin{center}
\section*{Introduction}
\end{center}

In recent years much work has been done in the study of infinite type hyperbolic manifolds, hyperbolizable manifolds with non-finitely generated fundamental group. For example lot of work has gone into studying the mapping class group of infinite type surfaces, for example \cite{AFP2017,Ba2016,PV2018}. Similarly, the first author has proven a hyperbolization result for a large class $\M^B$ of infinite type 3-manifolds, see \cite{C2017}. The class $\M^B$ is characterised by the fact that each $M\in\M^B$ has an exhaustion $\set{M_i}_{i\in\N}$ in which each $M_i$ is a compact, hyperbolizable 3-manifold with incompressible boundary and such that each $S\in\pi_0(\partial M_i)$ has genus at most $g=g(M)$. The class of hyperbolic 3-manifolds we will look at, denoted by $\M^{\S}$, are manifolds that need to admit exhaustions by $\pi_1$-injective sub-manifolds. Thus, we allow $M_i\subset M_{i+1}$ to have compressing disks in $M_i$, and we do not have any condition on the genus of the boundary components. However, we do need an embedding $\cup_{i\in\N} M_i\hookrightarrow \S$ and we will assume that $M\in\M^{\S}$ has no rank two cusps.

By work of Souto-Stover \cite{SS2013} and of Cremaschi-Souto \cite{CS2018} and Cremaschi \cite{C2018c,C20171} it is not hard to build hyperbolizable infinite type 3-manifolds that are homeomorphic to Cantor set complements in the 3-sphere $\S$. In particular, in \cite{CS2018}, the manifold of Example 2 can be extended to be a Cantor set complement showing, for example, how one can have a hyperbolizable Cantor set complement in $\S$ whose fundamental group is not residually finite.

The flexibility to build hyperbolizable Cantor set complements in $\S$ is reminiscent of the fact that most knots in $\S$ are hyperbolic. For example, of the $1,701,936$ knots with fewer than $16$ crossings, all but $32$ are hyperbolic, see \cite{HTW1998}. Moreover, Purcell-Souto \cite{PS2010} showed that if $M\hookrightarrow \S$ is a one-ended hyperbolic 3-manifold of finite type without parabolics, then $M$ is the geometric limit of hyperbolic knot complements. This shows how, under the geometric topology, hyperbolic knots are dense in the space of one ended-hyperbolic 3-manifolds admitting embeddings in $\S$.

The aim of the present work is to show a similar statement for hyperbolic 3-manifolds, not necessarily of finite type, admitting an embedding in $\S$. As approximating manifolds we will use Cantor sets complements:

\begin{customthm}{\ref{mainthm}} Let $M\cong \hyp\Gamma$ be a hyperbolic 3-manifold, not necessarily of finite type, without rank two cusps admitting an embedding $\iota:M\hookrightarrow \S$. Then, there exists a sequence of Cantor sets $\mathcal C_i\subset \S$, $i\in\N$, such that:
							\begin{itemize}
							\item[(i)] $N_i\eqdef \S\setminus \mathcal C_i$ is hyperbolic $N_i\cong\hyp{\Gamma_i}$;
							\item[(ii)] the $N_i$ converge geometrically to $M$.
							\end{itemize} 
\end{customthm}

As in \cite{PS2010}, one can obtain hyperbolic Cantor set complements with small eigenvalues of the Laplacian, arbitrarily large isometrically embedded balls, arbitrarily many short geodesics or surfaces with arbitrarily small principal curvatures.

\vspace{0.3cm}

\paragraph{Acknowledgements:} The first author would like to thank J.Souto for suggesting the problem.

\section{Background}

				\subsection{Notation and Conventions} All appearing 3-manifolds are assumed to be aspherical and orientable.  We use $\cong$ for homeomorphic. By $S\hookrightarrow M$, we denote an embedding of $S$ into $M$ while $S\imm M$ denotes an immersion. By $\Sigma_{g,n}$ we denote an orientable surface of genus $g$ with $n$ boundary components. We say that a manifold is closed if it is compact and without boundary. By $\pi_0(M)$ we denote the set of connected components of $M$, and unless otherwise stated we use $I=[0,1]$ to denote the closed unit interval.

Let $M$ be an open manifold, by an \emph{exhaustion} $\set{M_i}_{i\in\N}$ we mean a nested collection of compact sub-manifolds $M_i\subset\text{int}(M_{i+1})$ with $\cup_{i\in\N}M_i=M$. By \emph{gaps} of an exhaustion $\set{M_i}_{i\in\N}$ we mean the connected components of $\overline {M_i\setminus M_{i-1}}$. We will use $\hat\C$ to denote the Riemann sphere.
						
\subsection{Some 3-manifold topology} 

We now recall some facts and definitions about 3-manifold topology. For more details on the topology of 3-manifolds, some references are \cite{He1976,Ha2007,Ja1980}.

Let $M$ be an orientable 3-manifold, then $M$ is said to be \emph{irreducible} if every embedded sphere $\mathbb S^2$ bounds a 3-ball $\mathbb B^3$. Given a connected properly immersed surface $S\imm M$, we say it is \emph{$\pi_1$-injective} if the induced map on the fundamental groups is injective. Furthermore, if $S\hookrightarrow M$ is embedded and $\pi_1$-injective we say that the surface $S$ is \emph{incompressible} in $M$. By the Loop Theorem \cite{He1976,Ja1980} if $S\hookrightarrow M$ is a two-sided surface that is not incompressible, we have that there is an embedded disk $D\subset M$ such that $\partial D=D\cap S$ and $\partial D$ is non-trivial in $\pi_1(S)$. Such a disk is called a \emph{compressing} disk. 

An irreducible 3-manifold with boundary $(M,\partial M)$ is said to have \emph{incompressible} \emph{boundary} if every map of a disk: $(\mathbb D^2,\partial\mathbb D^2)\hookrightarrow (M,\partial M)$ is homotopic via maps of pairs into $\partial M$. Therefore, a manifold $(M,\partial M)$ has incompressible boundary if and only if each component $S$ of $\partial M$ is incompressible. We say that a 3-manifold $M$ is \emph{atoroidal} if any $\pi_1$-injective torus $T\subset M$ is homotopic into $\partial M$.

\bdefi
We say that a 3-manifold $M$ is \emph{hyperbolic}, or hyperbolizable, if $M\cong\hyp\Gamma$ for $\Gamma\subset\text{PSL}_2(\C)$ a discrete and torsion-free subgroup. The group $\Gamma$ is called \emph{Kleinian}.
\edefi 

In general hyperbolic 3-manifolds that are not closed are open. We will make use of the following convention:

\vspace{0.3cm}

\begin{centering} \emph{If we say that a compact 3-manifold is hyperbolic, we mean the interior and if $M$ is a finite type hyperbolic 3-manifold, we use $\overline M$ to mean its compact manifold closure.}
\end{centering}

\vspace{0.3cm}

The above convention makes sense since by the Geometrization Theorem \cite{Kap2001} and Tameness Theorem \cite{AG2004,CG2006}, any hyperbolic 3-manifold $M$ with finitely generated fundamental group is homeomorphic to the interior of a compact 3-manifold $\overline M$ such that $\overline M$ is irreducible, atoroidal and has infinite fundamental group.

Given a hyperbolic 3-manifold $M\cong\hyp\Gamma$ the \emph{convex core} $CC(M)\subset M$ is the smallest submanifold with convex boundary that whose inclusion induces a homotopy equivalence to $M$. We say that $M\cong\hyp \Gamma$ is \emph{convex co-compact} if $CC(M)$ is a compact submanifold and we say that $M$ is \emph{geometrically finite} if $CC(M)$ has finite volume. Some reference for hyperbolic 3-manifolds are: \cite{Ma2007,MT1998,Th1978,BP1992}.

	We now prove a couple of topological Lemmas.
	
			\blem\label{incboundary} Let $M$ be a compact 3-manifold with boundary and let $\P\subset \partial M$ be a collection of pairwise disjoint simple closed curves containing a pants decomposition of $\partial M$. Let $\gamma_i$, $1\leq i\leq n$, be the components of $\P$ and assume that every $\gamma_i$ is $\pi_1$-injective in $M$. For $0<g_i<\infty$, $1\leq i\leq n$, let $M'$ be the 3-manifold obtained by attaching thickenings of $\Sigma_{g_i,1}$ to $M$ by identifying regular neighbourhoods in $M$ of $\gamma_i$ and $\partial\Sigma_{g_i,1}$. Then, $\partial M'$ has incompressible boundary.
			\elem
			\bpf
			Without loss of generality, we can assume that $M$ has connected boundary. Let $(D,\partial D)$ be a compressing disk for $(M',\partial M')$. By an isotopy of $D$ we can assume that $D\pitchfork U$ for $U$ a regular neighbourhood of $\P$ in $\partial M$. 
			
If $U\cap D=\emp$ we have that $\partial D\subset \partial M\setminus \P$ hence $D$ is either in $M$ or in some $\Sigma_{g_i,1}$. If $D\subset M$, since $\P$ contains a pants decomposition, it means that $\partial D$ is isotopic into $\P$ giving us a contradiction with the fact that each component of $\P$ $\pi_1$-injects in $M$. If $D$ is contained in some $ \Sigma_{g_i,1}\times I$ we have that $\partial D\subset \Sigma_{g_i,1}\times\partial I$ but $\Sigma_{g_i,1}\times\partial I$ has no compressing disks in the $I$-bundle.		
			
Therefore, we have that $\mathcal A\eqdef D\cap U$ is a, non-empty, collection of essential arcs. Let $D'\subset D$ be an innermost disk with respect to the arc system $\mathcal A\subset D$. Then, $D'\cap U$ has only one component in $\partial D'$. Since $\P$ contains a pants decomposition, up to an isotopy of $D'$, we obtain a disk in either $M$ or $\Sigma_{g_i,1}\times I$ intersecting $U$ in an essential arc $\alpha$.

The disk $D'$ cannot be contained in $\Sigma_{g_i,1}\times I$ because every compressing disk intersects $\partial \Sigma_{g_i,1}\times I$ in at least two components. If $D'\subset M$ then $\partial D'$ is decomposed into two arcs $\alpha,\beta$ with $\alpha $ an essential arc in $U$ and $\beta$ an essential arc in $\partial M\setminus U$. However, since $\P$ contains a pants decomposition and $U$ is a thickening of $\P$, there cannot be such an essential $\beta$. \epf

Our last preliminary topological lemma
\blem\label{handlebodydisk}  Let $M$ be a compact 3-manifold with non-empty boundary $\partial M$ such that no component of $\partial M$ is a torus. Given, $\iota: M\hookrightarrow \S$ with handle-body complement $H$, we can find a pants decomposition $\P$ of $\partial M$ such that $\P$ is a disk-system for $H$ and is $\pi_1$-injective in $M$. \elem
\bpf
 Let $\mathcal D$ be a disk system\footnote{Such a disk system always exists and is even possible given a disk system $\mathcal D$ to surger it, by band sums, to get a new disk system $\mathcal D'$ that has no separating component.} for $H$ such that no disk $D\in\mathcal D$ is separating in $H$. We now need to show that the loops $\iota^{-1}(\partial \mathcal D)$ are essential in $M$. If not, by the loop Theorem if $\gamma$ in $\iota^{-1}(\partial \mathcal D)$ is not $\pi_1$-injective in $M$, then it bounds a disk $D'$. Let $D$ be the disk of $\mathcal D$ corresponding to $\gamma$. Then $S=D\cup_\gamma D'$ is an embedded $2$-sphere in $\S$, and so it is separating. However, since each $D$ is non-separating in $H\subset\S$ we get a contradiction. \epf

\brem
In the setup of Lemma \ref{incboundary} and \ref{handlebodydisk} we can take a disk system so that no pair is separating in $\partial H$ and so that the manifold $M'$ of Lemma \ref{incboundary} has incompressible boundary and the JSJ decomposition of $M'$ is given by the thickened surfaces we attach.
\erem

\subsection{Combination Theorems}

For the reader's convenience, we now recall some Theorems dealing with glueings of Kleinian groups, i.e. hyperbolic 3-manifolds.

\bthm[\cite{Kap2001}, 4.97]\label{pingpong}
Let $G_1$, $G_2$ be Kleinian groups with fundamental domains $D_1$, $D_2$ in $\hat \C$ such that: $\hat\C\setminus D_2\subset\text{int}(D_1)$ and $\hat\C\setminus D_1\subset\text{int}(D_2)$. Then, the group $G$ generated by $G_1$ and $G_2$ is Kleinian and isomorphic to $G_1*G_2$. Moreover, $D\eqdef D_1\cap D_2$ is a fundamental domain for $G$ on $\hat \C$.
\ethm

\bdefi
Let $\Gamma\subset PSL_2(\C)$ be a Kleinian group. Given a subgroup $H\subset \Gamma$ we say that $B\subset\widehat \C$ is \emph{precisely invariant} under $H$ if $H(B)\subset B$ and for all $\gamma\in\Gamma\setminus H$ we have that $\gamma(B)\cap B=\emp$.
\edefi

\bthm[\cite{Kap2001}, 4.104]\label{cyclamalg1}
Let $G_1$, $G_2$ be a pair of Kleinian groups such that $G_1\cap G_2=H$, where $H$ is a cyclic subgroup. Let $D_j$ be fundamental domains for the actions of $G_j$ on $\hat \C$, $j=1,2$. Let $B_1$, $B_2$ be open disks in $\hat \C$ such that $J\eqdef \overline B_1\cap \overline B_2=\partial B_1=\partial B_2$ is a topological circle. Suppose the following:
\begin{itemize}
\item $B_j$ is invariant under $H$ in $G_j$, $j=1,2$;
\item $D_j'\eqdef D_j\cap G_j(B_j)\subset B_j$, $j=1,2$;
\item $D_1'\cap D_2$ and $D_1\cap D_2'$ have non-empty interiors.
\end{itemize}
Then, the subgroup $G\subset \text{Isom}(\mathbb H^3)$ generated by $G_1$, $G_2$ is Kleinian and isomorphic to $G_1*_HG_2$. If $G_1$, $G_2$ are geometrically finite, then $G$ is also geometrically finite. The quotient $\Omega(G)/G$ is naturally conformally equivalent to:
$$\Omega(G_1\setminus G_1(B_1))/G_1\cup_L\Omega(G_2\setminus G_2(B_2))/G_2$$
where the glueing is along $L=[J\cap \Omega(H)]/H$. Any parabolic element in $G$ is either conjugate to $G_1$ or to $G_2$ or conjugate to an element commuting with a parabolic element of $H$.
\ethm

Similarly:

\bthm[\cite{Kap2001},4.105]\label{cyclamalg2}
Let $G_0$ be Kleinian and $H_1$, $H_2$ a pair of cyclic subgroups. Let $D_0$ be a fundamental domain for the actions of $G_0$ on $\hat\C$. Let $B_1$, $B_2$ be open disks in $\hat \C$ and $A\in\text{Isom}(\mathbb H^3)$ be a M\"obius transformation such that $AH_1A^{-1}=H_2$. This conjugation induces an isomorphism $\phi:H_1\rar H_2$. Suppose the following:
\begin{itemize}
\item $B_j$ is precisely invariant under $H_j$ in $G_0$, $j=1,2$;
\item $A(B_1)\cap B_2=\emp$ and $A(\partial B_1)\cap\partial B_2=J$ is a topological circle;
\item $gB_1\cap B_2=\emp$ for all $g\in G_0$;
\item $D_0\cap (\hat\C\setminus G_0(B_1\cup B_2))$ has non-empty interior.
\end{itemize}
Then, the subgroup $G\subset \text{Isom}(\mathbb H^3)$ generated by $G_0$, $A$ is Kleinian and isomorphic to the HNN-extension $G_0*_{\phi:H_1\rar H_2}$ of $G_0$ via $\phi$. If $G_0$ is geometrically finite, then $G$ is also geometrically finite. The quotient $\Omega(G)/G$ is naturally conformally equivalent to:
$$ \sim/ [\Omega(G_0)\setminus G_0(B_1\cup B_2)]/G_0$$
where the identification is such that $[J\cap \Omega(H_2)]/H_2$ is identified with $[A^{-1}(J)\cap \Omega(H_1)]/H_1$ via the projection of $A$. Any parabolic element in $G$ is either conjugate to $G_0$ or conjugate to an element commuting with a parabolic element of $H_j$, $j=1,2$.
\ethm

\brem[Parabolic amalgamation]\label{paramalg} Let $z$ be a parabolic fixed point for the action of a Kleinian group $\Gamma$ corresponding to a 3-manifold $M$. By the Universal Horoball Theorem \cite[3.3.4]{Ma2007}, we can always find an embedded horoball $H$ in $\Omega(\Gamma)$. Therefore, by using the universal horoball, it is easy to glue Kleinian groups $\Gamma_1$ and $\Gamma_2$ along a common parabolic group $\langle \alpha\rangle$.\erem

\section{Reduction to the Convex Co-Compact case}

We start by recalling a useful Lemma about converging sequences of geometric limits.
	\blem\label{geomsub} If $M$ is the geometric limit of $\set{M_i}_{i\in\N}$ and each $M_i$ is the geometric limit of $\set{N_i^n}_{n\in\N}$, then $M$ is the geometric limit of a sub-sequence $\set{N^n_{a_n}}_{n\in\N}$. 
	\elem
\bpf
Consider the diagram:
\be\xymatrix{ (M_i,p_i)\ar[r]^i & (M,p)\\ 
(N_i^n,q^n_i)\ar[u]_n\ar@{.>}[ur]}\ee
By geometric convergence in $i$ we have that $\forall R>0:\exists i_R$ such that $\forall i\geq i_R$ we have embeddings:
\be f_i:(B_R(p),p)\hookrightarrow (M_i,p_i)\qquad f_i\;\;(1+\epsilon_i)\text{-bilipschitz}\quad \epsilon_i\rar 0\ee
and similar statements for $(N_i^n,q_i^n)$ and $(M_i,p_i)$.

For each $i$ we have that $f_i(B_R(p))\subset B_{R+\epsilon_i}(p_i)$ thus we have $(1+\zeta_{i,n})$-bilipschitz embeddings $g_i^n:B_{R+\epsilon_i}(p_i)\rar (N^n_i,q^n_i)$. Therefore, the embeddings:
\be g^n_i\circ f_i:B_R(p)\rar (N_i^n,q^n_i)\ee
are $(1+\epsilon_i)(1+\zeta_{i,n})$-bilipschitz. Thus, we can find a geometrically convergent sub-sequence.
\epf

We now reduce the general case to the convex co-compact case. 
\bdefi We say that a 3-manifold $M$ is in $\M^{\S}$ if $M\hookrightarrow \S$ is hyperbolic without rank two cusps: $M\cong\hyp\Gamma$, $\Gamma\subgroup PSL_2(\C)$ and $M$ is either of finite type, i.e. $\pi_1(M)$ is finitely generated or $M=\cup_{i\in\N} M_i$ (as in the above sense) in which $\pi_1(M_i)\hookrightarrow \pi_1(M)$. The last  condition is equivalent to, up to sub-sequence, $\pi_1(M_i)\hookrightarrow \pi_1(M_{i+1})$.
\edefi

\blem\label{approx1} Let $M\cong \hyp\Gamma$ be a hyperbolic 3-manifold, not necessarily of finite type and with $\Gamma$ not abelian, without rank two cusps, admitting an embedding $\iota:\overline M\hookrightarrow \S$. If $M$ admits an exhaustion by $\pi_1$-injective compact sub-manifolds then, there is a sequence of finite type hyperbolic 3-manifolds with no parabolics $(M_i,p_i)$ 3-manifolds with embeddings $f_i:(\overline M_i,p_i)\hookrightarrow\S$ such that $(M_i,p_i)\rar (M,p)$ geometrically.
\elem
\bpf  Let $N_i$'s be the $\pi_1$-injective sub-manifolds giving us an exhaustion of $M$, and let $\Gamma_i\subset \Gamma$ be the corresponding Kleinian groups. Without loss of generality, we can assume that $N_i\not\simeq N_{i+1}$ so that $\Gamma_i\neq\Gamma_{i+1}$. Then, $\Gamma_i\subsetneq \Gamma_{i+1}$ and $\cup_{i\in\N}\Gamma_i=\Gamma$. Then, we obtain the required sequence by: 
$$(M_i,p_i)\eqdef\left(\hyp{\Gamma_i},[0]\right).$$

Since the $N_i$ are $\pi_1$-injective in $M$, they lift homeomorphically to the covers $\pi_i:M_i\rar M$. By Tameness \cite{AG2004,CG2006} we have that $M_i\setminus N_i$ are product regions and so $M_i\cong \text{int}(N_i)$.  Hence, the $M_i$ also embed in $\S$, concluding the proof.\epf

\bprop \label{approx2}Let $M\cong \hyp\Gamma$ be a hyperbolic 3-manifold in $\M^{\S}$. Then, there is a sequence of convex co-compact hyperbolic 3-manifolds $(M_i,p_i)$ with embeddings $f_i:(\overline M_i,p_i)\hookrightarrow\S$ such that $(M_i,p_i)\rar (M,p)$ geometrically.
\eprop
\bpf We first deal with the case $\Gamma$ is abelian, hence of finite type. Any such Kleinian group can be geometrically approximated by a classical Schottky group on two generators and we are done. 

Let $(M_i,p_i)$ be the sequence from Lemma \ref{approx1}. Since each $M_i$ has no $\mathbb Z^2\in\pi_1(M_i)$ by the Strong Density Theorem \cite[1.4]{BS06}, there is a collection of convex co-compact manifolds $N_n^i\in AH(M_i)$ converging strongly to $M_i$, moreover without loss of generality, by geometric convergence, we can assume that for all $n:$ $N_n^i\cong M_i$. By Lemma \ref{geomsub}, we have a sub-sequence $N_{n_i}^i$ that converges geometrically to $M$. Moreover, since each $N_{n_i}^i$ is homeomorphic to $M_i$, they admit embeddings $f_i:\overline N_{n_i}^i\rar\S$. \epf
\brem
 The previous proposition is the only place in the paper in which we actually need the exhaustion and the fact that we have no rank two cusps.
\erem

\subsection{General Proof assuming convex co-compact approximation}

We now assume the following Theorem, which we will prove in the next sections. The main step will be a glueing argument that is done in Section 3.

\begin{customthm}{\ref{convexcocompact}} Let $M\cong \hyp\Gamma$ be a convex co-compact hyperbolic 3-manifold admitting an embedding $\iota:\overline M\hookrightarrow \S$. Then, there exists a sequence of Cantor sets $\mathcal C_i\subset \S$, $i\in\N$, such that:
							\begin{itemize}
							\item[(i)] $N_i\eqdef \S\setminus \mathcal C_i$ is hyperbolic $N_i\cong\hyp{\Gamma_i}$;
							\item[(ii)] the $N_i$ converge geometrically to $M$.
							\end{itemize} 
\end{customthm}
and prove:

\bthm\label{mainthm} Let $M\cong \hyp\Gamma$ be a hyperbolic 3-manifold and let $M\in\M^{\S}$. Then, there exists a sequence of Cantor sets $\mathcal C_i\subset \S$, $i\in\N$, such that:
							\begin{itemize}
							\item[(i)] $N_i\eqdef \S\setminus \mathcal C_i$ is hyperbolic $N_i\cong\hyp{\Gamma_i}$;
							\item[(ii)] the $N_i$ converge geometrically to $M$.
							\end{itemize} 
\ethm
\bpf By Proposition \ref{approx2}, we have a sequence of convex co-compact manifolds $\overline M_i\hookrightarrow \S$ that converge geometrically to $M$. By Theorem \ref{convexcocompact}, each $M_i$ is approximated by Cantor set complements; hence, by Lemma \ref{geomsub} $M$ is approximated, geometrically, by Cantor set complements.
\epf

\section{Gluing argument}

In this section, we will show how given $M\subset \S$ convex co-compact such that $M^C$, the complement of $M$ in $\S$, is a collection of handlebodies $H$ we can extend the metric of $M$ to a new 3-manifold $M'$ such that $M\subsetneq M'\subset \S$ and $H'\eqdef (M')^C$ is a collection of handlebodies such that $H'\subsetneq H$. Moreover, each component of $H$ contains at least two components of $H'$, and for $h\in \pi_0(H)$ and $h'\in\pi_0(H')$ we have: $\diam(h')\leq \frac 12 \diam (h)$. By iterating this argument, we will build our hyperbolic Cantor set complements. The aim of this section is to show our main glueing argument:

\begin{customprop}{ \ref{glueing1}}Let $M$ be a convex co-compact hyperbolic manifold  with the property that $\P\subset \partial M$ is a $\pi_1$-injective collection of pairwise disjoint simple closed curves. Let $m\eqdef \abs{\P}$ and let $L\in [0,\infty)$. Then, there exists $\set{g_i}_{i=1}^m$ with $1\leq g_i<\infty$ such that we can extend the hyperbolic metric of $M$ to a convex co-compact manifold: $$M_L\eqdef M \cup_\P \coprod_{i=1}^m \Sigma_{g_i,1} \times I$$ with the property that:
\begin{enumerate}
\item in $\Sigma_{g_1,1}\times I$, the geodesic corresponding to $\P_i$ has a collar of width at least $L$;
\item if $\P$ contains a pants decomposition, then $M_L$ has incompressible boundary.
\end{enumerate}
\end{customprop}

Before showing Proposition \ref{glueing1}, we show that given a compact convex co-compact manifold $M$ embedding in $\S$, we can assume, up to geometric limit, that it has handle-body complement.

\blem\label{hyphandlebodies} Let $\iota: M\hookrightarrow \S$ be a compact convex co-compact hyperbolic manifold. Then, by adding a collection of $1$-handles $H$ to $M$, we have an embedding $\iota':M\cup_\partial H\hookrightarrow \S$, extending the metric, such that $\overline{\S\setminus \iota'(M\cup_\partial H)}$ is a collection of handlebodies and $M\cup_\partial H$ is convex co-compact.
\elem
\bpf If $\overline{\iota (M)^C}$ is a collection of handlebodies, there is nothing to do. Otherwise, let $N\subset \overline{\iota (M)^C}$ be a non-handlebody component. Let $C=\mathcal H_g\cup Q$ be a minimal genus Heegaard splitting of $N$, where $\mathcal H_g$ is a genus $g$ handlebody, and $Q$ is a collection of $2$-handles. Attaching a 2-handle $P$ to $H_g$ is equivalent to attaching a $1$-handle $P'$ to $\overline{\iota (M)}$. Thus, we get that by attaching all $1$-handles to $\iota (M)$ we can make $N$ a handlebody component. Therefore, there is a collection of  $1$-handles $H$ and an embedding $\iota': M\cup H\hookrightarrow \S$ such that $\overline{\S\setminus \iota'(M\cup_\partial H)}$ is a collection of handlebodies.

We now need to show that we can realise the above topological construction while extending the given hyperbolic metric on $M\cong\hyp \Gamma$. This essentially follows from Ping-Pong Lemma (Theorem \ref{pingpong}). There are two cases depending on whether the $1$-handle $P$ is attached to one or two boundary components of $M$. We will indicate by $S_1$ and $S_2$ these two boundary components.

Assume that $S_1\neq S_2$. Let $D_1$ be a fundamental domain for the action of $\Gamma$ on $\hat\C$. Since $\Gamma$ is convex co-compact $\Gamma.D_1$ has full measure and let $F_1\eqdef\hat \C\setminus D_1$. Pick two points $x_1$ and $x_2$ in $\text{int}(D_1)\cap \tilde S_1$ and $\text{int}(D_1)\cap \tilde S_2$ respectively and let $h_\lambda\in\text{ Isom}^+(\textbf H^3)$, $\lambda\in(0,\infty)$, be the loxodromic element with fixed points $x_1$ and $x_2$ and translation length $\lambda$. Let $D_2(\lambda)$ be the fundamental domain of $\langle h_\lambda\rangle$ and $F_2\eqdef \hat \C\setminus D_2$. Since as $\lambda\rar\infty$: 
\be D_2(\lambda)\xrightarrow{\text{Hausdorff}}\hat\C\setminus \set{x_1,x_2}\qquad F_2(\lambda)\Orar{\text{Hausdorff}} \set{x_1,x_2} \ee
 we get that there is $\lambda\in(0,\infty)$ such that:
\be D_2(\lambda)\supset F_1 \qquad D_1\supset F_2(\lambda) \ee
Then, by Theorem \ref{pingpong} $\Gamma'\eqdef\langle \Gamma, h_\lambda\rangle$ is discrete, isomorphic to $\Gamma*h_\lambda$ and $\hyp{\Gamma'}$ has the required topological type.

If $S_1=S_2$ let $D_1$ and $F_1$ as before and pick $x\neq y$ to be points in $D_1\cap \tilde S_1$. Then, by the same reasoning as before, we can find $h_\lambda$ such that $\Gamma'\eqdef\langle \Gamma, h_\lambda\rangle$ is discrete, isomorphic to $\Gamma*h_\lambda$ and $\hyp{\Gamma'}$ has the required topological type.
\epf

We now define: 
\bdefi Let $N$ be a geometrically finite 3-manifold, we say that the convex core of $N$ is homeomorphic to $\Sigma_{g,k,n}\times I$ if $CC(N)$ has $n$ rank $1$ cusps, $k$ funnels and there is a type-preserving homeomorphism $f:N\xrightarrow{\cong} \Sigma_{g,k,n}\times I$.
\edefi

The next Lemma constructs a handlebody piece that will be attached to $M$ via cyclic amalgamation, Theorem \ref{cyclamalg1}, along a peripheral loxodromic $\gamma$. This particular construction produces a rank-1 cusp that we will have to deal with later. The loxodromic element $\gamma$ and $\gamma$-invariant disk $B\subset\partial_\infty \Ha^3$ in the statement will be obtained from $M$ by taking an incompressible curve in $\partial M$ and lifting a collar around it.

\blem\label{handlemodel1} Given $\gamma\in PSL_2(\C)$ loxodromic element and a closed $\gamma$-invariant disk $B\subset \partial_\infty \Ha^3$ then there is a Schottky group extension of $\langle\gamma\rangle$, $\Gamma = \Gamma_B$, such that
\begin{enumerate}
    \item\label{limitset} The limit set of $\Gamma$ is included in $B$.
    \item\label{handletopology} The convex core of $\hyp{\Gamma}$ is homeomorphic to $\Sigma_{g,1,1}\times [0,1]$, where the boundary component of $\Sigma_{g,1,1}$ corresponds to $\gamma$ and the puncture to a rank-1 cusp.
\end{enumerate}
Moreover, such group $\Gamma$ can be taken so that $\gamma$ has a collar larger than any given constant.
\elem

\bpf
Take $B'\subset B$, a smaller region delimited by two $\gamma$-invariant smooth arcs $\rho_1,\rho_2$ joining the fixed points of $\gamma$. Furthermore, select a third $\gamma$-invariant smooth path $\rho\subset B'$ so that $A=B'/\langle\gamma\rangle$ is an annulus with boundary $\rho_1/\langle\gamma\rangle\cup \rho_2/\langle\gamma\rangle$ and $\pi_1$ representative embedded curve $\rho/\langle\gamma\rangle$. 

\begin{center}\begin{figure}[hb!]
						\def\svgwidth{160pt}
						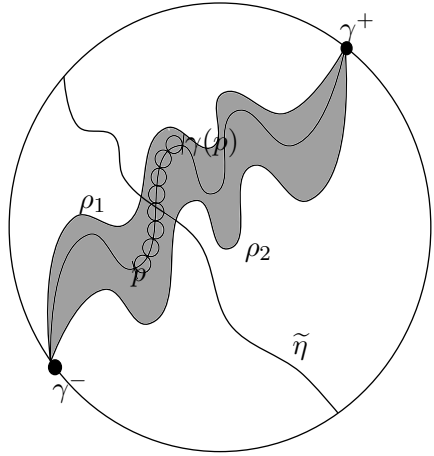
						\caption{The disks $\Delta_i$ in the domain of discontinuity containing a lift of $\gamma$. The shaded region is the ball $B'$}\label{circlepacking}
						\end{figure}
\end{center}

We would to find $F\subset B'$ so that $F$ is a fundamental region of $A$ and $\rho\cap F = \rho_F$ is connected. To do this, denote by $\gamma_\pm$ the fixed points of $\gamma$. Consider a closed path $\eta$ in the annulus $(\partial_\infty\Ha^3\setminus\lbrace\gamma_\pm\rbrace)/\langle\gamma\rangle$, such that $\eta$ intersects each one of $\rho/\langle\gamma\rangle, \rho_1/\langle\gamma\rangle, \rho_2/\langle\gamma\rangle$ exactly once (and the annulus $A=B'/\langle\gamma\rangle$ in a connected segment). Define $F_0$ as the lift in $\partial_\infty\Ha^3$ of the complement of $\eta$ in $(\partial_\infty\Ha^3\setminus\lbrace\gamma_\pm\rbrace)/\langle\gamma\rangle$. This makes $F_0$ a disjoint union of connected components, and the closure of any of these components is a fundamental domain for $(\partial_\infty\Ha^3\setminus\lbrace\gamma_\pm\rbrace)/\langle\gamma\rangle$. Then one can verify that $F$ can be obtained by $F:=F_0\cap B$. Cover $\rho_F$ by closed disks $\lbrace\Delta_i\rbrace_{-N\leq i\leq N}$ in $B'$, see Figure \ref{circlepacking}, such that:

\begin{enumerate}
    \item $\Delta_i, \Delta_{i+1}$ are tangent for all $0\leq i\leq 4N-2$, $\lbrace p_{i+1}\rbrace = \Delta_i \cap \Delta_{i+1}$.
    \item $\Delta_i\cap\Delta_j=\emptyset$ for $|i-j|\geq 2$.
    \item $\Delta_{4N-1}=\gamma(\Delta_{0})$
\end{enumerate}

Iterate by powers of $\gamma$ to obtain, $\lbrace\Delta_i\rbrace_{i\in\mathbb{Z}}$, a covering of $\rho$ by disk in $B$ such that:

\begin{enumerate}
    \item $\Delta_i, \Delta_{i+1}$ are tangent for all $i\in\mathbb{Z}$, $\lbrace p_{i+1}\rbrace = \Delta_i \cap \Delta_{i+1}$.
    \item $\Delta_i\cap\Delta_j=\emptyset$ for $|i-j|\geq 2$.
    \item $\Delta_{i+4N}=\gamma(\Delta_{i})$
\end{enumerate}

Select $f_i$ a Mobius map that sends the triple $(\partial_\infty \Ha^3 \setminus (\mathring\Delta_i) ,p_{i} , p_{i+1})$ to the triple $(\Delta_{i+2}, p_{i+3},p_{i+2})$, so that $\gamma^{-1}\circ f_i\circ \gamma = f_{i+4N}$ (make a priori such a selection). Furthermore, denote by $a_i=f_{4i}$, $b_i=f_{4i+1}$.  Let $\Gamma$ be the group generated by $a_0, b_0,\ldots, a_{N-1}, b_{N-1},\gamma$ (also generated by $\langle\lbrace a_i, b_i \rbrace_{i\in\mathbb{Z}},\gamma\rangle$). Then by modifying the proof of Theorem \ref{pingpong} we can prove that $\Gamma$ is a Kleinian group freely generated by $a_0, b_0,\ldots, a_{N-1}, b_{N-1},\gamma$. Indeed, for $a_i, b_i$ take fundamental domains as the complement of the appropriate disks $\Delta_i$, and take $F_0$ as the fundamental domain for $\gamma$. Taking any two of these fundamental domains (and denoting them by $D_1, D_2$) we have that $D_1\supset(\mathbb{C}\setminus D_2), D_2\supset(\mathbb{C}\setminus D_1)$ rather than $\text{int}(D_1)\supset(\mathbb{C}\setminus D_2), \text{int}(D_2)\supset(\mathbb{C}\setminus D_1)$, the latter as in Theorem \ref{pingpong}. This is because of the tangencies we consider. Nevertheless, if $D$ denotes the intersection of all fundamental domains and $w$ is a nontrivial word generated by $a_0, b_0,\ldots, a_{N-1}, b_{N-1},\gamma$, it follows that for any $z\in\text{int}(D)$ we have $w(z)\notin D$. This implies that $\Gamma$ is freely generated by $a_0, b_0,\ldots, a_{N-1}, b_{N-1},\gamma$ and that $D$ is a fundamental domain for $\Gamma$ in $\partial_\infty\mathbb{H}^3 = \mathbb{S}^2$.

Remains to show that $\Gamma$ is discrete. Take $x\in \mathbb{H}^3$ in the complement of all the half-spaces bounded by the $\Delta_i$'s, intersected with a fundamental domain of $\gamma$ bounded by $F_0$ (which is the complement of two topological half-spaces). Take into consideration that all half-spaces can be taken mutually disjoint. Then assume that there is a sequence $\lbrace g_k\rbrace\subset \Gamma$ so that $g_k(x)\rightarrow x$. For any $g_k\neq id$ we have that $g_k(x)$ belongs to one of the discarded half-spaces. Hence $g_k=id$ for $k$ sufficiently large, and from which we know that $\Gamma$ is a Kleinian group. And since $D$ is a fundamental domain for $\Gamma$ in $\partial_\infty\mathbb{H}^3 = \mathbb{S}^2$, it follows that the limit set of $\Gamma$ is contained in $B$. This is because the complement of $\langle\gamma\rangle D$ is contained in $B'$.

Note that all the points of tangencies $\lbrace p_i\rbrace$ are identified with one another in the quotient by $\Gamma$, where the element $\alpha = \gamma[a_{N-1},b_{N-1}],\ldots[a_0,b_0]$ fixes $p_{0}$. Moreover, $\alpha$ preserves the direction tangential to the disks meeting at $p_{0}$. In order to make $\alpha$ parabolic, we can make choices so that $D\alpha_{p_{0}}$ has norm $1$ with respect to the standard $\mathbb{S}^2$ metric. Take the loxodromic element $c_\lambda$ with real translation and fixed points $p_2, p_3$, so that the derivatives of $c_\lambda$ at $p_2, p_3$ are $\lambda^{-1}, \lambda$, respectively. We can choose then $c_\lambda\circ a_0$ instead of $a_0$. The new choice $c_\lambda\circ a_0$ satisfies the same conditions as $a_0$ and introduces a factor $\lambda$ twice while applying chain rule for $D\alpha_{p_{0}}$ (once for $c_\lambda$ at $p_3$ and once for $c_\lambda^{-1}$ at $p_2$). Then by taking the appropriate value for $\lambda$ we make $\alpha$ parabolic. We claim then that such $\Gamma$ is geometrically finite with convex core homeomorphic to $\Sigma_{g,1,1}\times [0,1]$, where the boundary component of $\Sigma_{g,1,1}$ corresponds to $\gamma$ and the puncture to the rank-1 cusp generated by $\alpha$. Indeed, we can select a smooth metric in $D/\Gamma$ so that $p_{0}$ is a hyperbolic cusp. By taking the Epstein envelope surface \cite{Epstein84} of a sufficiently small multiple of the selected metric, we obtain a finite volume core with convex boundary. Then $\Gamma$ is geometrically finite. Finally, the boundary of the core can be easily seen as $\Sigma_{2g,0,2}$, where $\gamma$ is a separating curve that divides the quotient into components homeomorphic to $\Sigma_{g,1,1}$. From here we can see that the convex core of $\hyp{\Gamma_B}$ is homeomorphic to $\Sigma_{g,1,1}\times [0,1]$, where $\gamma$ corresponds to the boundary component of $\Sigma_{g,1,1}$ and the puncture to a rank-1 cusp.

As a final remark, observe that the collar around $\gamma$ gets bigger as we take the region $B$ and the disks $\Delta_i$ smaller. \epf

We now start the first step of our main glueing construction:

\blem\label{glueing0} Let $M$ be a convex co-compact hyperbolic manifold with the property that $\P\subset \partial M$ is a $\pi_1$-injective collection of disjoint non-homotopic curves. Let $n\eqdef \abs{\P}$ and let $L\in [0,\infty)$. Then, there exists $\set{g_i}_{i=1}^n$ with $1\leq g_i<\infty$ such that we can extend the hyperbolic metric of $M$ to a geometrically finite manifold: $$M_L'\eqdef M \cup_\P \coprod_{i=1}^n \Sigma_{g_i,1,1} \times I$$ with the property that:
\begin{enumerate}
\item $\Sigma_{g_1,1,1}\times I$ has a rank one cusp corresponding to a boundary component of $\Sigma_{g_1,1,1}$, and the other boundary is glued to a component $\P_i$ of $\P$;
\item in $\Sigma_{g_1,1,1}\times I$, the geodesic corresponding to $\P_i$ has a collar of width at least $L$.
\end{enumerate}
\elem
\bpf
Since each element $\P_i$ in $\P$ is $\pi_1$-injective, then it has a loxodromic element $\gamma_i\in\pi_1(M)$, and a $\gamma_i$-invariant disk $B_i$ in the domain of discontinuity of $M$. By Lemma \ref{handlemodel1} there exist Schottky group extensions $\Gamma_{B_i}$ with limit set in $B_i$ and collars around $\gamma_i$ as large as we desire. Then by Theorem \ref{cyclamalg1} the manifold $M'=M'(L)$ with fundamental group generated by $\langle\pi_1(M),\Gamma_{B_1},\ldots\Gamma_{B_n}\rangle$ has the desired properties, provided that the groups $\lbrace \Gamma_{B_i}\rbrace$ from Lemma \ref{handlemodel1} have all collars bigger than $L$ around the geodesics that each of them is extending. \epf

\begin{center}\begin{figure}[h!]
						\def\svgwidth{200pt}
						%% Creator: Inkscape 1.0 (4035a4f, 2020-05-01), www.inkscape.org
%% Accompanies image file 'cantoset2.pdf' (pdf, eps, ps)
\begingroup%
  \makeatletter%
  \providecommand\color[2][]{%
    \errmessage{(Inkscape) Color is used for the text in Inkscape, but the package 'color.sty' is not loaded}%
    \renewcommand\color[2][]{}%
  }%
  \providecommand\transparent[1]{%
    \errmessage{(Inkscape) Transparency is used (non-zero) for the text in Inkscape, but the package 'transparent.sty' is not loaded}%
    \renewcommand\transparent[1]{}%
  }%
  \providecommand\rotatebox[2]{#2}%
  \newcommand*\fsize{\dimexpr\f@size pt\relax}%
  \newcommand*\lineheight[1]{\fontsize{\fsize}{#1\fsize}\selectfont}%
  \ifx\svgwidth\undefined%
    \setlength{\unitlength}{455.8644068bp}%
    \ifx\svgscale\undefined%
      \relax%
    \else%
      \setlength{\unitlength}{\unitlength * \real{\svgscale}}%
    \fi%
  \else%
    \setlength{\unitlength}{\svgwidth}%
  \fi%
  \global\let\svgwidth\undefined%
  \global\let\svgscale\undefined%
  \makeatother%
  \begin{picture}(1,0.53340335)%
    \lineheight{1}%
    \setlength\tabcolsep{0pt}%
    \put(0,0){\includegraphics[width=\unitlength,page=1]{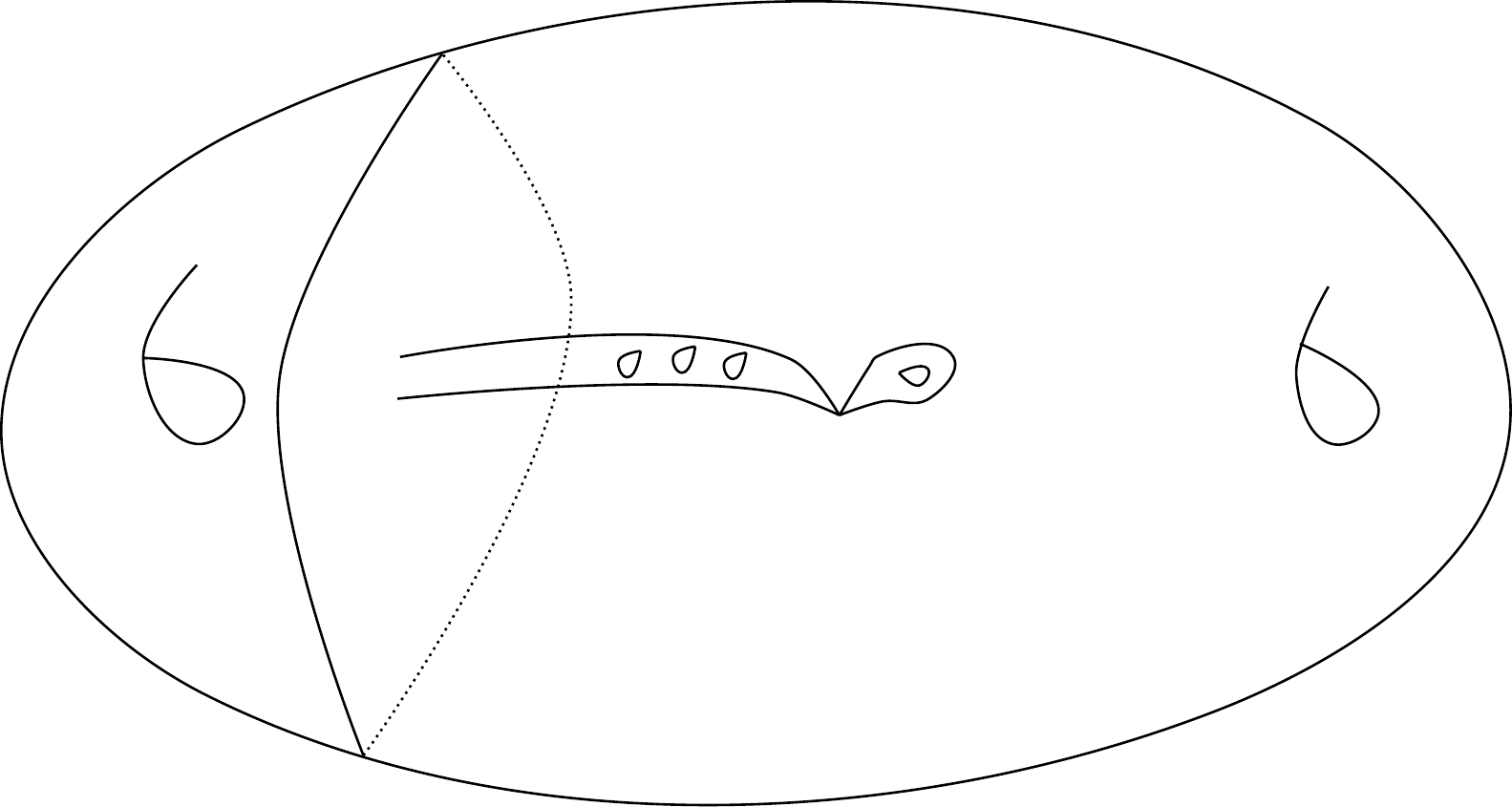}}%
    \put(0.21873244,0.0478058){\color[rgb]{0,0,0}\makebox(0,0)[lt]{\lineheight{1.25}\smash{\begin{tabular}[t]{l}$\gamma_i$\end{tabular}}}}%
  \end{picture}%
\endgroup%

						\caption{Partial stage in which we glued a punctured torus to a $\Sigma_{3,1,1}\times I$ along one $\gamma_i$. The rank $1$ cusp $C_i$ and accidental parabolic correspond to the node.}\label{pic1}
						\end{figure}
\end{center}

We can now prove our main glueing step, where we will deal with the parabolics:

\bprop \label{glueing1}Let $M$ be a convex co-compact hyperbolic manifold  with the property that $\P\subset \partial M$ is a $\pi_1$-injective collection of pairwise disjoint simple closed curves. Let $m\eqdef \abs{\P}$ and let $L\in [0,\infty)$. Then, there exists $\set{g_i}_{i=1}^m$ with $1\leq g_i<\infty$ such that we can extend the hyperbolic metric of $M$ to a convex co-compact manifold: $$M_L\eqdef M \cup_\P \coprod_{i=1}^m \Sigma_{g_i,1} \times I$$ with the property that:
\begin{enumerate}
\item in $\Sigma_{g_i,1}\times I$ the geodesic corresponding to $\P_i$ has a collar of width at least $L$;
\item If $\P$ contains a pants decomposition, then $M_L$ has incompressible boundary.
\end{enumerate}
\eprop
\bpf Start with the manifold $M_L'$ coming from Lemma \ref{glueing0} and let $C_i$ be the rank $1$ cusps corresponding to the $\Sigma_{g_i,1,1}\times I$ attached to $\gamma_i$. By applying Klein-Maskit combination (Theorem \ref{cyclamalg1}) to universal horoballs to each rank $1$ cusps, we attach a $\Sigma_{1,1}\times I$ manifold. This gives us a new manifold:
$$M_L''\eqdef M \cup_\P \coprod_{i=1}^m \Sigma_{g_i+1,1} \times I$$
in which the $ \Sigma_{g_i+1,1} \times I$ have an accidental parabolic $\delta_i$ corresponding to the remaining rank 1 cusp $C_i$ coming from the Klein-Maskit combination, see Figure \ref{pic1}.

Note that if $\P$ contains a pants decomposition, then, by Lemma \ref{incboundary}, the manifold $M_L''$ has incompressible boundary.

For each rank-1 cusp, we can find invariant tangent disk at the corresponding fixed point, by cyclic amalgamation, see Remark \ref{paramalg}, we can glue each rank-1 cusp onto itself to produce a geometrically finite manifold with rank-2 cusps. Each cusp has an embedded cylinder towards each of the two boundary components where it appears as an accidental parabolic.

Thus, we get manifolds: 
$$\widehat M_L= M \cup_\P \coprod_{i=1}^m \left(\Sigma_{g_i+1,1} \times I\setminus \delta_i\times\set {1/2}\right)$$
still extending the metric on $M$.

By Thurston's Dehn Filling Theorem \cite{Th1978,BP1992,Co1996}, we have $N\in \N$ such that for all $n>N$ the manifolds $\widehat M_L^n$ obtained from $\widehat M_L$ by doing $\frac 1 n$-Dehn Filling on every rank two cusp, see \cite{KT1990}, are convex co-compact. Moreover, by taking a larger $N$, if necessary, we can also assume that: 
$$\widehat M_L^n\cong M \cup_\P \coprod_{i=1}^m \left(\Sigma_{g_i+1,1}\times I\right)$$
where the homeomorphisms $\phi_n$ restrict to the identity on $M$ and are induced by $\tau_{\gamma_i}^n$, the $n$-th Dehn twist along $\delta_i$, on $\Sigma_{g_i+1,1}\times I$. Hence, for all $L$ and $n$, the manifolds $\widehat M_L^n$ are convex co-compact and have incompressible boundary by Lemma \ref{incboundary}.

Finally, we have that:
$$\widehat M_L^n\overset{geom}{\underset{n\rar\infty}{\longrightarrow}} \widehat M_L.$$
Thus, by definition of geometric convergence, by taking $n$ large enough and some $L'>L$, we can assume that in $M_L\eqdef \widehat M_{L'}^n$ all the geodesics corresponding to $\P$ have a collar of width at least $L$. Hence, the manifold $M_L$ satisfies all the requirements of the proposition completing the proof. \epf

\bcor\label{glueing2} Let $M$ be a convex co-compact hyperbolic manifold with the property that $\P\subset \partial M$ is a $\pi_1$-injective collection of pairwise disjoint simple closed curves. Let $m\eqdef \abs{\P}$, $p\in CC(M)$, $R>0$ and $n\in\N$ there exists $L=L(p,R,n)$ and: 
$$f: N_R( CC(M))\hookrightarrow M \cup_\P \coprod_{i=1}^m \left(\Sigma_{g_i+1,1}\times I\right)$$
such that $f$ is $(1+\frac 1 n)$-bi-Lipschitz.  \ecor
\bpf Pick $\set{L_n}_{n\in\N}\subset\R^+$ such that $L_n\nearrow \infty$. Build the manifolds $M_n\eqdef M_{L_n}$ as in Proposition \ref{glueing1}. It is easy to see that for any $p\in CC(M)$, by property (1) of Proposition \ref{glueing1}, the sequence:
$$(M_n,p)\overset{geom}{\longrightarrow}(M,p)$$
giving us the desired result.
\epf

We can now prove our iteration step. One of the main takeaways is that we can choose our embeddings so that the diameter of the complement decays to 0 as we iterate the process, which is necessary to obtain a Cantor set complement.

\bprop \label{cantoriteraton} Let $\iota: M\hookrightarrow \S$ be an embedding of a compact irreducible manifold whose complement is a collection of handlebodies $H_i$, $1\leq i\leq n$. Then, by attaching a finite collection $\Sigma_{g_h^i,1}\times I$, $1\leq h \leq n_i$ to a collection of disks $D_i$ on $\partial H_i$, containing a disk system for $H_i$, we obtain a new embedding:
$$\iota': M\bigcup_{i=1}^n \cup_{h=1}^{n_i}\Sigma_{g_h^i,1}\times I\hookrightarrow \S$$
extending $\iota$ such that $\iota'(\cup_{h=1}^{n_i}\Sigma_{g_h,1})\subset H_i$ and $\overline{H_i\setminus{\iota'(\cup_{h=1}^{n_i}\Sigma_{g_h^i,1})}}$ is a collection $J_1^i,\dotsc, J_{m_i}^i$ of handlebodies with $m_i\geq 2$ and $\diam(J^i_{m_j})\leq\frac 1 2 \diam(H_1)$. Moreover, if $\text{int}(M)\cong \hyp\Gamma$ is convex co-compact given $L>0$, we can extend the hyperbolic metric to $M\bigcup_{i=1}^n \cup_{h=1}^{n_i}\Sigma_{g_h^i,1}$ so that each attaching region has a collar of width at least $L$.
\eprop
\bpf

Let $\Gamma$ be the hyperbolic structure on $M$. It suffices to prove the statement for each handlebody component $H_i$, for the sake of notation, we will just refer to it as $H$. Let $\mathcal D$ be the disk system coming from Lemma \ref{handlebodydisk}.

\begin{center}\begin{figure}[h!]
						\def\svgwidth{250pt}
						%% Creator: Inkscape 1.0 (4035a4f, 2020-05-01), www.inkscape.org
%% Accompanies image file 'cantoset3.pdf' (pdf, eps, ps)
\begingroup%
  \makeatletter%
  \providecommand\color[2][]{%
    \errmessage{(Inkscape) Color is used for the text in Inkscape, but the package 'color.sty' is not loaded}%
    \renewcommand\color[2][]{}%
  }%
  \providecommand\transparent[1]{%
    \errmessage{(Inkscape) Transparency is used (non-zero) for the text in Inkscape, but the package 'transparent.sty' is not loaded}%
    \renewcommand\transparent[1]{}%
  }%
  \providecommand\rotatebox[2]{#2}%
  \newcommand*\fsize{\dimexpr\f@size pt\relax}%
  \newcommand*\lineheight[1]{\fontsize{\fsize}{#1\fsize}\selectfont}%
  \ifx\svgwidth\undefined%
    \setlength{\unitlength}{415.1184082bp}%
    \ifx\svgscale\undefined%
      \relax%
    \else%
      \setlength{\unitlength}{\unitlength * \real{\svgscale}}%
    \fi%
  \else%
    \setlength{\unitlength}{\svgwidth}%
  \fi%
  \global\let\svgwidth\undefined%
  \global\let\svgscale\undefined%
  \makeatother%
  \begin{picture}(1,0.38716198)%
    \lineheight{1}%
    \setlength\tabcolsep{0pt}%
    \put(0,0){\includegraphics[width=\unitlength,page=1]{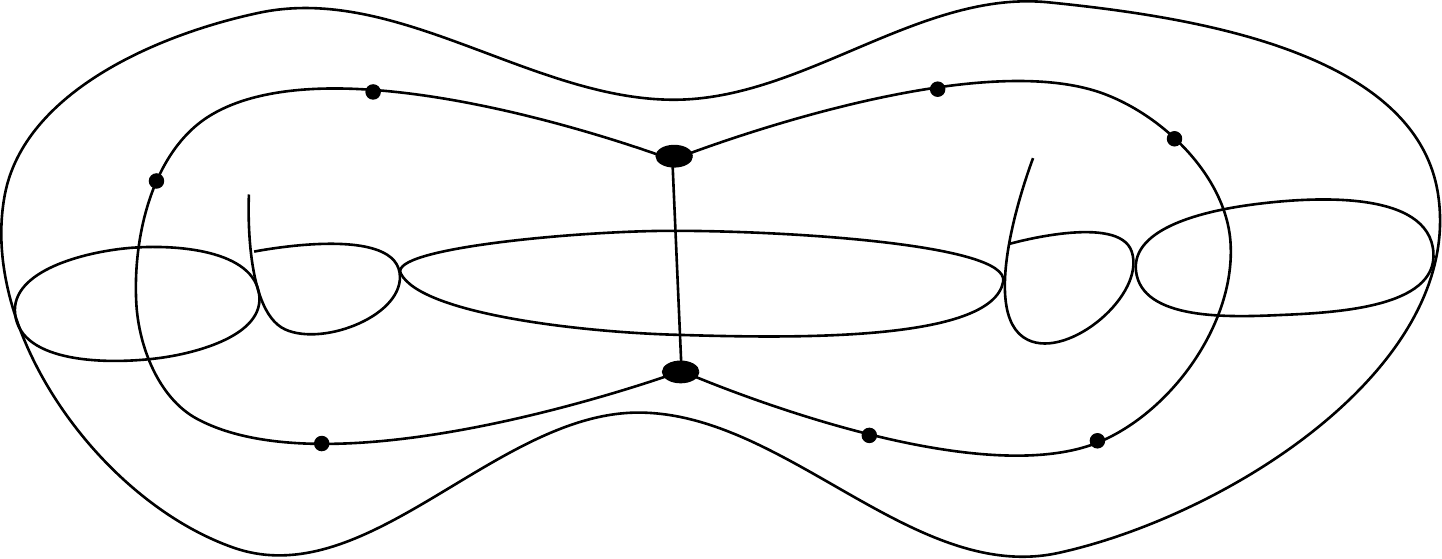}}%
    \put(0.15829931,0.99258993){\color[rgb]{0,0,0}\makebox(0,0)[lt]{\lineheight{1.25}\smash{\begin{tabular}[t]{l}--\end{tabular}}}}%
    \put(0,0){\includegraphics[width=\unitlength,page=2]{cantoset3.pdf}}%
    \put(0.63815687,0.08573611){\color[rgb]{0,0,0}\makebox(0,0)[lt]{\lineheight{1.25}\smash{\begin{tabular}[t]{l}$\ell_j$\end{tabular}}}}%
  \end{picture}%
\endgroup%

						\caption{Nerve subdivision. The shaded ball $B_j$ is the thickening of an $\ell_j$ section whose diameter is less than $\frac14\diam(H)$.}\label{img:diameter}
						\end{figure}
\end{center}

Take a nerve on the handlebody $H$ so that in each ball component of $\overline{H\setminus N_\epsilon(\mathcal D)}$ we have a trivalent vertex. By using copies of disks in $\mathcal D$, we subdivide the nerve into sections $\ell_1,\ldots,\ell_\kappa$ so that each ball component $B_m$, $1\leq m\leq \kappa$, has diameter less than $\frac14\diam(H)$, see figure \ref{img:diameter}. 

This gives us a collection of disks $\mathcal D'\subset H$ containing a pants decomposition of $\partial H$. Moreover, each component of $\mathcal D'$ $\pi_1$-injects in $M$. Then, by applying Corollary \ref{glueing2} to $(\mathcal D',\partial\mathcal D')$ we obtain a hyperbolic 3-manifold $M\cup_{h=1}^{n}\Sigma_{g_h,1}$ extending $\Gamma$. 

We now construct the nested family of handlebodies obtained by successively attaching handles to the curves homotopic to $\partial \mathcal{D}'$. We do this so that each section from $\ell_1,\ldots,\ell_\kappa$ appear inside a handlebody.  To each disk $D\in\pi_0(\mathcal D')$, we attach $g_h$ 2-handles by drilling them from the adjacent 3-ball in an unknotted way so that they complement is a handlebody.

Each handlebody $J_1,\dotsc, J_\kappa$ is a thickening of an element of $\ell_1,\ldots,\ell_\kappa$ with some handles attached or drilled in. Moreover, we can do it so that the resulting handle is still close to the corresponding element of $\ell_1,\ldots,\ell_\kappa$, and more importantly so that each complementary region's diameter is less than $\frac12\diam(H)$. Since $\kappa\geq 3g(H)-3>2$, we complete the proof.	\epf

\brem  Note that we can make the resulting manifold of Proposition \ref{cantoriteraton} boundary incompressible by selecting a $\pi_1$-injective pants decomposition during the last iteration of the handle attaching.
\erem

\section{Proof for Convex Co-Compact}

Before proving the main result, we prove the following key Proposition:

\bprop\label{cantorlimit}
Let $(M,p)$ be a convex co-compact hyperbolic 3-manifold admitting an embedding $\iota:(\overline M,p)\hookrightarrow \S$ with complement given by a collection of handlebodies $\mathcal H$. Given $R>0$, there exists a Cantor set $C_R\subset \mathcal H$ such that $\S\setminus C_R$ is hyperbolizable and $B_R(p)\subset \S\setminus C_R$ is $1+e(R)$ bilipschitz to the $R$-ball around $p$ in $M$. Moreover, $e(R)\rar 0$ as $R\rar \infty$.
\eprop 
\bpf Pick $L>R$ and apply Proposition \ref{cantoriteraton} to $\iota:M\hookrightarrow\S$ to obtain a new manifold $N_1^L$ so that all new topology is at distance $L>R$ from $CC(M)$. We then reiterate this construction using the same $L$. We thus obtain a collection of convex co-compact hyperbolic 3-manifolds $N_n^L$ admitting nested embedding $\overline N_n^L\subset\overline N_{n+1}^L$ whose complement in $\S$ is a collection of handlebodies $H_n$ and whose direct limit $N_\infty^L$ is homeomorphic to the complement of a sub-set $K$ of $\S$.

\textbf{Claim 1:} The set $K$ is a Cantor set so that $N_\infty^L\cong \S\setminus \mathcal C_R$.

\bpfc To show that $K$ is a Cantor set, we need to show that it is a compact, perfect, totally disconnected metric space. Let $C\eqdef \diam (H_1)$. By construction, it is easy to see that $K=\cap_{n\in\N} H_n$ where each $H_n$ is a collection of handlebodies in which each component of $H_n$ contains at least two components of $H_{n+1}$. Moreover, by Proposition \ref{cantoriteraton}, we have that for $H$ a component of $H_n$: $\diam H\leq 2^{-n} C$ so that $K$ is a collection of points. Since each component of $H_n$ contains at least two components of $H_{n+1}$ we see that $K$ is also totally disconnected. Thus, being a closed sub-set of a compact metrisable space, it is compact and metrisable as well. The fact of it being perfect is also a straightforward consequence of the nesting construction. \epfc

\textbf{Claim 2:} The $B_R(p)\subset \S\setminus C_R$ is $1+e(L)$ bi-lipschitz to the $R$-ball around $p$ in $M$ and $e(L)\rar 0$ as $L\rar \infty$.

\bpfc This follows from Proposition \ref{glueing1}. \epfc

If $R\rar \infty$, so does $L$, and the last claim of the Proposition is proven. \epf

We now finish the proof of the main result:

\bthm\label{convexcocompact}  Let $M\cong \hyp\Gamma$ be a convex co-compact hyperbolic 3-manifold admitting an embedding $\iota:\overline M\hookrightarrow \S$. Then, there exists a sequence of Cantor sets $\mathcal C_i\subset \S$, $i\in\N$, such that:
							\begin{itemize}
							\item[(i)]$N_i\eqdef \S\setminus \mathcal C_i$ is hyperbolic $N_i\cong\hyp{\Gamma_i}$;
							\item[(ii)] the $N_i$ converge geometrically to $M$.
							\end{itemize} 
\ethm	\bpf By Lemma \ref{hyphandlebodies}, we can assume that we have $M_i\rar M$ geometrically with embeddings $\iota_i:\overline M_i \rar \S$ such that $\iota_i(\overline M_i)^C$ are handlebodies for every $i$. Then, by Lemma \ref{geomsub}, it suffices to prove the Theorem for such an $M_i$.

Thus, let $M$ be a convex co-compact hyperbolic 3-manifold with an embedding $\iota:\overline M \rar \S$ that has for complement a collection of handlebodies $\mathcal H=\set{H_1,\dotsc, H_n}$.

Choose any strictly increasing sequence $R_n$. By applying Proposition \ref{cantorlimit} to $(M,p,R_n)$ we obtain a sequence of Cantor set complements $(\S\setminus \mathcal C_n,p, R_n)$ that geometrically converge to $M$, concluding the proof.\epf

Since, in particular, $\mathbb H^3\hookrightarrow \S$ we have Cantor sets complements $N_n\eqdef \S\setminus \mathcal C_n$ and points $p\in \mathbb H^3$ and $p_n\in N_n$ such that:
$$(N_n,p_n)\rar (\mathbb H^3,p)$$
geometrically.  Thus, the balls of radius $B_R(p)\subset \mathbb H^3$ can be $(1+\epsilon_n)$-isometrically embedded in $N_n$. In particular, this means that for large enough $n$ the set of points of distance, say, $\frac R2$ from $p_n$ is simply connected and so $inj_{p_n}(N_n)\geq \frac R2$. Since $R$ was arbitrary, we obtain:
\bcor For all $R>0$, there exists a Cantor sets $\mathcal C\subset\S$ such that $\S\setminus\mathcal C$ is hyperbolic and there is a point $p\in \S\setminus\mathcal C$ with injectivity radius at least $R$.
\ecor 

However, we do not necessarily know what the shape of the corresponding Cantor set is. Moreover, as in \cite{PS2010}, one can obtain hyperbolic Cantor set complements with small eigenvalues of the Laplacian, arbitrarily many short geodesics or surfaces with arbitrarily small principal curvatures.

\nocite{BP1992,CEM2006,CM2006,Ha2002,He1976,Sh1975,Sc1972,SY2013,Th1978,Ja1980,MT1998,Ca1993,MT1998,So2006,Th1982}

\thispagestyle{empty}
{\small
\markboth{References}{References}
\bibliographystyle{abbrv}
\bibliography{mybib}{}
}
	\bigskip
	\small
	{\noindent Department of Mathematics, University of Southern California, Kaprelian Hall.

\noindent 3620 S. Vermont Ave, Los Angeles, CA 90089-2532.

email: \texttt{cremasch@usc.edu}

\vspace{0.1cm}

\noindent Department of Mathematics,  Yale University.

 \noindent 12 Hillhouse, New Haven, CT 06511
 
 email: \texttt{franco.vargaspallete@yale.edu}
}
\end{document}